\documentclass[10pt,a4paper,twoside,leqno]{article}
\usepackage{amsfonts,amssymb}
\usepackage{eufrak}
\usepackage{amscd}

\newtheorem{defn}{Definition}

\vspace{0.5cm}

 4
\font\ebf=cmbx8
\font\erm=cmr8

\usepackage{graphicx}

\begin{document}

\thispagestyle{empty}

\noindent {\bf Cauchy $\hat q_{\psi}$-identity and $\hat q_{\psi}$
-Fermat matrix via $\hat q_{\psi}$-muting variables of $\hat
q_{\psi}$-Extended Finite Operator Calculus}

\vspace{0.7cm} {\it Andrzej K. Kwa\'sniewski}

\vspace{0.7cm}

{\erm Higher School of Mathematics and Applied Informatics}

\noindent{\erm  Kamienna 17, PL-15-021 Bia\l ystok , Poland}

\vspace{0.5cm}

\noindent {\ebf Summary}

{\small New $\psi$-labeled family of Cauchy identities is found
and Fermat matrix notion -now with operator entries - is also
$\psi$-extended including standard Cauchy  and $q$-Cauchy or new
Fibonomial-Cauchy cases due to the use of $\hat
q_{\psi}$-commuting variables introduced and promoted as
representative for Extended Finite Operator Calculus by the author
few years ago.

\vspace{0.7cm}

\section{I. Towards The Non-commuting }

At first we make a $\psi$- remark on the notation. $\psi$ denotes
an extension of $\langle\frac{1}{n!}\rangle_{n\geq 0}$ sequence to
quite arbitrary one (the so called - admissible) and the specific
choices are for example: Fibonomialy-extended ($\langle F_n
\rangle$ - Fibonacci sequence )
$\langle\frac{1}{F_n!}\rangle_{n\geq 0}$   or just "the usual"
$\langle\frac{1}{n!}\rangle_{n\geq 0}$ or Gauss $q$-extended
$\langle\frac{1}{n_q!}\rangle_{n\geq 0}$ admissible sequences of
extended umbral operator calculus. We get used to write these $q$
- Gauss and other extensions in mnemonic convenient  upside down
notation [3-6]
\begin{equation}\label{eq1}
\psi_n\equiv n_\psi , x_{\psi}\equiv \psi(x)\equiv\psi_x ,
 n_\psi!=n_\psi(n-1)_\psi!, n>0 ,
\end{equation}
\begin{equation}\label{eq2}
x_{\psi}^{\underline{k}}=x_{\psi}(x-1)_\psi(x-2)_{\psi}...(x-k+1)_{\psi}
\end{equation}
\begin{equation}\label{eq3}
x_{\psi}(x-1)_{\psi}...(x-k+1)_{\psi}= \psi(x)
\psi(x-1)...\psi(x-k-1) .
\end{equation}
You may consult [3-6]  and references therein for further
development and use of this notation  "$q$-commuting variables" -
included.
    The idea to use "$q$-commuting variables" goes back at least to Cigler (1979) [1]
( see formula  (7) , (11) in [1] ) and also Kirchenhofer - see [2]
for further systematic development ) We shall take here notation
from [3-6] and the results from [3-6] - for granted- in view of an
easy access via ArXiv to the source papers . For other respective
references to Cigler, Kirchenhofer and Others see: [3-6]. The task
is to invent what to replace with the question mark in the Cauchy
identity type formula $(4)$ below

\begin{equation}\label{eq4}
\sum_{k\geq 0}  (?)\Big({r \atop k}\Big )_{\psi}\Big({s
\atop{j-k}}\Big )_{\psi}= \Big({{r+s} \atop j}\Big )_{\psi}   ,
\Big({n \atop k}\Big )_{\psi}=0 ,{ for } k<0 .
\end{equation}

to get true new extended identities . The simple minded way to do
it is closed because for

$$ (x+_{\psi} 1)^r(x+_{\psi} 1)^s = (\sum_{k\geq 0}\Big({r \atop k}\Big
)_{\psi}x^k)(\sum_{l\geq 0}\Big({s \atop l}\Big )_{\psi}x^l) $$

$$ E(\partial _\psi )x^{r+s} \equiv (x+_\psi 1)^{r+s}\equiv \sum_{j\geq 0}\Big({{r+s} \atop j}\Big
)_{\psi}x^j   $$

the regular arrive at the  Cauchy identity fails  as $$(x
+{\psi}y)^r(x +{\psi}y)^s \neq (x +{\psi}y)^{r+s} .$$

 For example: $$(x +_F y)^2 = x^2 + F_2xy + y^2, (x +_F y)^4 = x^4 + F_4 x^3y +
F_4F_3x^2y^2 + F_4 xy^3 + y^4.$$
$$(x +_F y)^5 = x^5 + F_5 x^4y + F_5F_4 x^3y^2 +...+ F_5xy^4 + y^5.$$

$$(x +_F y)(x +_F y)^4 \neq (x +F y)^5  . $$
 where $\{{F_n}\}_{n>0}$  denotes the  Fibonacci sequence  (see subsequent pages).
What then would help perhaps [3-6] is to replace commuting
variables $xy = yx$ by non-commuting ones as done by Cigler (see:
also Kirchenhofer) in Umbral Calculus domain. The idea to use
"$q$-commuting variables" goes back at least to Cigler (1979)
[1,2] (see formula (7), (11) in [1]) and also Kirchenhofer - see
[2] for further systematic development). In [2] Kirchenhofer
equivalently defined  the polynomial sequence of $q$ -binomial
type by [5,4]

\begin{equation}\label{eq5}
p_{n} \left( {A + B} \right) \equiv \sum\limits_{k \ge 0} \Big({r
\atop k}\Big )_q p_{k} \left( A \right)p_{n - k} \left( B \right)
{ where }  [B,A]_q \equiv BA-qAB =0.
\end{equation}

 $A$ and $B$ might be interpreted here as
co-ordinates on quantum $q$-plane (see [7] Chapter 4). For example
$A = \hat {x}$ and $B = y \hat {Q}$ where $\hat {Q}\varphi \left(
{x} \right) = \varphi \left( {qx} \right)$ , (more on that [5,4]).

In the $q$-case (see: Proposition 4.2.3 in [7]- with $yx = qxy$)
we have
\begin{equation}\label{eq6}
\sum_{k\geq 0}  q^{(r-k)(j-k)}\Big({r \atop k}\Big )_q\Big({s
\atop{j-k}}\Big )_q= \Big({{r+s} \atop j}\Big )_q
\end{equation}

hence  from this  Cauchy  $q$-identity here we define immediately
$(s = i , s = j)$  the symmetric $q$-Pascal (or $q$-Fermat) matrix
elements via the following  easy to find out formula:

\begin{equation}\label{eq6}
\sum_{k\geq 0}  q^{(r-k)(j-k)}\Big({i \atop k}\Big )_q
\Big({j\atop k}\Big )_q = \Big({{i+j} \atop j}\Big )_q .
\end{equation}

For arbitrary admissible sequences labeling the elements of the
giant family of Extended Finite Operator Calculus -
\textbf{calculi } of Rota , Roman and Others   the present author
introduced in  [4,5] appropriate (also for our purpose here)
notions -to be used next. Note: \textbf{Calculi} means "stones"
(or "pebbles" or "counters") in Latin. And really indeed - some
choices of admissible $\psi$-sequences are impressive
combinatorial.

\section{II. Examples of choices and $\hat q_\psi $-binomial symbol}

\vspace{2mm}

Let since now on  $\psi= \{\psi_n(q)\}_{n\geq 0}$ as in [3-6].
With the Gauss  choice  $\psi_n(q)=[n_q!]^{-1}$ where $
n_q=\frac{1-q^n}{1-q}          $
 and $n_q!=n_q(n-1)_q! , 1_q!=0_q!=1$ we may interpret $q$-binomial coefficients
 in a standard way, namely: $q$-Gaussian coefficient $\Big({n \atop k}\Big )_q $
 denote number of $k$-dimensional subspaces in $n-th$ dimensional space over
Galois field  $GF(q)$ [8,9] i.e. we are dealing with lattice of
subspaces. For  $q=1$ we arrive at the  lattice of subsets and the
binomial coefficient standard interpretation. In [8] the
combinatorial interpretation  was proposed  also for  Fibonomial
coefficients
$$
\left( \begin{array}{c} n\\k\end{array}
\right)_{F}=\frac{F_{n}!}{F_{k}!F_{n-k}!}\equiv
\frac{n_{F}^{\underline{k}}}{k_{F}!},\quad n_{F}\equiv F_{n}\neq
0, $$

\noindent where we make an analogy driven [6,4,3] identifications
$(n>0)$:
$$
n_{F}!\equiv n_{F}(n-1)_{F}(n-2)_{F}(n-3)_{F}\ldots 2_{F}1_{F};$$
$$0_{F}!=1;\quad n_{F}^{\underline{k}}=n_{F}(n-1)_{F}\ldots (n-k+1)_{F}. $$

 (Here:  $\psi_n(q)=[F_n!]^{-1}$ ). In [8,9]  a partial ordered set was defined in such a way that
 the Fibonomial coefficients count the number of specific finite
"birth-self-similar"  sub-posets of this  infinite non-tree poset
naturally related to the Fibonacci tree of rabbits growth process.
For fascinating "weighted choices" - Konvalina combinatorial
interpretations see [10,11]. In order to proceed we take from
[4,5]  only this what we need now.(Vector spaces are over the
field of zero characteristics). An so let us define the main
notion of this note [4,5].

\vspace{2mm}

\begin{defn}
Let $\left\{ {p_{n}}  \right\}_{n \ge 0} $ be the $\partial
_{\psi}$-basic polynomial sequence of the $\partial _{\psi}
$-delta operator $Q\left( {\partial _{\psi} }  \right) = Q$. Then
the $\hat {q}_{\psi ,Q} $-operator is a linear map\\ $\hat
{q}_{\psi ,Q} :P \to P; \quad \hat {q}_{\psi ,Q} p_{n} =
\frac{{\left( {n + 1} \right)_{\psi}  - 1}}{{n_{\psi} }
}p_{n},\quad n \ge 0$.
\end{defn}

\vspace{2mm}

We call this useful  $\hat {q}_{\psi ,Q} $ operator the $\hat
{q}_{\psi ,Q}$\textit{-mutator operator}.\\

\vspace{2mm}

 \textbf{Note:} For $Q =id\;$ $Q\left( {\partial
_{\psi} } \right) =
\partial _{\psi}  $ the natural notation is $\hat {q}_{\psi ,id}
 \equiv \hat {q}_{\psi}  $. For $Q =id$ and $\psi _{n} \left( {q}
\right) = \frac{{1}}{{R\left( {q^{n}} \right)!}}$ and $R\left( {x}
\right) = \frac{{1 - x}}{{1 - q}} \quad \hat {q}_{\psi ,Q}  \equiv
\hat {q}_{R,id}
 \equiv \hat {q}_{R} \equiv \hat {q}_{q,id}  \equiv \hat
{q}_{q}  \equiv \hat {q}$ and $\hat q x^{n}$ = $q^{n}x^{n}$.

\vspace{2mm}

\begin{defn}
Let $A$ and $B$ be linear operators acting on $P$; $A:P \to P$,
$B: P \to P$. Then $AB - \hat {q}_{\psi ,Q}BA \equiv  [A,B]_{\hat
{q}_{\psi ,Q}}  $ is called $\hat {q}_{\psi ,Q} $-mutator of $A$
and $B$ operators.
\end{defn}
Consider then the following special case $\hat {q}_{\psi} $ of
linear on $P= F[x]$ { $\hat {q}_{\psi ,Q}$ -mutator operator } :

$\hat {q}_{\psi} :P \to P; \quad \hat {q}_{\psi} x^n =
\frac{{\left( {n + 1} \right)_{\psi}  - 1}}{{n_{\psi} } }x^n
;\quad n \ge 0$.

Consider also $\hat {q}_{\psi}$ -muting variables $yx=\hat
q_{\psi} xy$  . Introduce also a -binomial symbol:

\vspace{2mm}

\begin{defn}
We define $\hat{q}_{\psi}$-binomial symbol i.e. $\hat
{q}_{\psi}$-Gaussian coefficients as follows: \\
$\Big({n \atop k}\Big)_{\hat {q}_{\psi}}=
\frac{n_{\hat{q}_{\psi}}!}{k_{\hat
q_{\psi}}!(n-k)_{\hat{q}_{\psi}}!}=\Big({n \atop {n-k}}\Big
)_{\hat {q}_{\psi}}$ { where } $n_{\hat {q}_{\psi}}!=n_{\hat
{q}_{\psi}}(n-1)_{\hat {q}_{\psi}}! , 1_{\hat {q}_{\psi}}!=0_{\hat
{q}_{\psi}}!=1$ and $n_{\hat {q}_{\psi}}=\frac{1-{\hat
{q}_{\psi}}^n}{1-{\hat {q}_{\psi}}}$ for $n>0$.
\end{defn}

\vspace{2mm}

\textbf{Challenge 1} Are we facing possibility of
$\hat{q}_{\psi}$-quantum "groups" investment alike  [7] ?

\section{III.  $\psi$-sequences labeled family of Cauchy identities}
Equipped with the above  we discover immediately the being looked
for Cauchy  $\psi$-identity formula and the $\psi$- extended
Fermat matrix (including  $q$- Fermat and Fibonomial $F$-Fermat
matrix cases). Namely let us observe that the following is true.

\vspace{2mm}

\textbf{Observation 1}
$$
p_{n} \left( {A + B} \right) \equiv \sum\limits_{k \ge 0} {\left(
{{\begin{array}{*{20}c}
 {n} \hfill \\
 {k} \hfill \\
\end{array}} } \right)} _{\hat{q}_{\psi}} p_{k} \left( {A} \right)p_{n - k} \left( {B}
\right) $$    where $$   [B,A]_{\hat{q}_{\psi}} \equiv
BA-\hat{q}_{\psi}AB = 0
$$
and in particular
$$
\left({x + y}\right)^n = \sum\limits_{k \ge 0} {\left(
{{\begin{array}{*{20}c}
{n} \hfill \\
{k} \hfill \\
\end{array}} } \right)} _{\hat{q}_{\psi}} x^k y^{n - k}
$$
where $$
[y,x]_{\hat{q}_{\psi}}=0.
$$

\vspace{2mm}

\textbf{Challenge 2}  Are we (compare with [5]) facing the
possibility of systematic representation of a general umbral
calculus  in Rota-like operator form  [3-6] with help of
$\hat{q}_{\psi}$-quantum "plane" variables in place of
"$q$-commuting variables" employed by Cigler (1979) [1] and
Kirchenhofer - to do their splendid efficient job - this time in

$\hat{q}_{\psi}$-extended binomial enumeration (this refers to
fundamental binomial enumeration formulation of umbra in [12])? -
Are we - ...?

In the $\hat{q}_{\psi}$-case of the $\hat{q}_{{\psi},Q}$ from
[5,6] general extended umbral theory in Rota-like finite operator
form we have

\vspace{2mm}

\textbf{Observation 2}
\begin{equation}\label{}
\sum_{k\geq 0}  \hat q_{\psi}^{(r-k)(j-k)}\Big({r \atop k}\Big )_{
\hat q_{\psi}}\Big({s \atop{j-k}}\Big )_{ \hat q_{\psi}}=
\Big({{r+s} \atop j}\Big )_{ \hat q_{\psi}}
\end{equation}

\vspace{2mm}

From Observation 2  i.e.  from the  Cauchy   $\hat q_{\psi}$ -
identity  we infer the following  $\hat q_{\psi}$ -formula for
matrix elements of the symmetric  $\hat q_{\psi}$- Pascal (or
$\hat q_{\psi}$-Fermat) matrix elements

\vspace{2mm}

\begin{equation}\label{}
\sum_{k\geq 0} {\hat q_{\psi}}^{(r-k)(j-k)}\Big({i \atop k}\Big
)_{\hat q_{\psi}}\Big({j\atop k}\Big )_{\hat q_{\psi}} =
\Big({{i+j} \atop j}\Big )_{\hat q_{\psi}} .
\end{equation}

\vspace{2mm}

For the first most recent applications see  [13]. For $q$-Pascal
matrix see [14].

In analogy to the standard case [15-17,14]  we shall call the
matrices (compare with [14,13]) -   with operator valued matrix
elements

$$
x^{i-j}\Big({i \atop j}\Big )_{\hat q_{\psi}}
$$
and
$$
\Big({{i+j} \atop j}\Big )_{\hat q_{\psi}}
$$
the  $\hat q_{\psi}$- Pascal P[x]  and  $\hat q_{\psi}$ -Fermat
F[1]  matrices - correspondingly.

\vspace{2mm}

The result from [15-17] tempt to be $\hat q_{\psi} $-extended.



\end{document}